\documentclass{amsart}
\usepackage{amssymb}
\usepackage{amsmath}
\usepackage{amsfonts}
\usepackage{hyphenat}
\usepackage[dvips]{graphicx}
\usepackage{float}

\newtheorem{theorem}{Theorem}[section]
\theoremstyle{plain}

\newtheorem*{definition}{Definition}

\newtheorem{lemma}[theorem]{Lemma}

\newtheorem{proposition}[theorem]{Proposition}
\newtheorem{remark}[theorem]{Remark}

\newtheorem*{main}{Main Theorem}
\numberwithin{equation}{section}
\newenvironment{proofnosquare}[1][Proof]{\smallskip\noindent\textsl{#1.} }{}
\newenvironment{Theorem Ref}[1][Theorem]{\medskip\noindent\textbf {#1}} {\medskip}
\begin{document}
\title{Strong S-equivalence of Ordered Links}
\author{Carol Gwosdz Gee}
\date{September 20, 2004}

\begin{abstract}
Recently Swatee Naik and Theodore Stanford proved that two
S-equivalent knots are related by a finite sequence of doubled-delta
moves on their knot diagrams. We show that classical S-equivalence
is not sufficient to extend their result to ordered links. We define
a new algebraic relation on Seifert matrices, called Strong
S-equivalence, and prove that two oriented, ordered links $L$ and
$L^{\prime }$ are related by a sequence of doubled-delta moves if
and only if they are Strongly S-equivalent. We also show that this
is equivalent to the fact that $L^{\prime }$ can be obtained from
$L$ through a sequence of Y-clasper surgeries, where each clasper
leaf has total linking number zero with $L$.
\end{abstract}

\maketitle

\section{Introduction}

A fundamental problem in knot theory is to classify knots and links by type,
and several algebraic and geometric invariants have been developed with this
goal in mind. \ Recently Swatee Naik and Ted Stanford have shown an
equivalence between S-equivalence of knots, a purely algebraic invariant,
and the doubled-delta move on knot diagrams, a purely geometric relation. \
They prove that two knots are related by a sequence of doubled-delta moves
if and only if the knots are S-equivalent \cite{ns}. \ We prove that the
analogous connection for links between S-equivalence and the doubled-delta
move, posed as a question by Stavros Garoufalidis, is false. \

We define a new invariant of ordered links called \emph{Strong S-equivalence}
that is in many ways better suited for links than the classical definition
of S-equivalence. \ With this new definition, we are able to prove a theorem
analogous to Naik and Stanford's result for links. \ Our main theorem also
ties these results to the emerging subject of clasper surgery, which has
strong connections to the field of \textquotedblleft quantum
topology,\textquotedblright\ an area of study that encompasses the Jones
polynomials, Vassiliev invariants, and the Kontsevitch integral.

\begin{main}
Consider two oriented, ordered $m$-component links $L_{0}$ and
$L_{1}$.  The following four statements are equivalent:

i. $ L_{1}$ can be obtained from $L_{0}$ through a sequence of
doubled-delta moves.

ii. $\ L_{0}$ and $L_{1}$ are related by a sequence of $Y$-clasper
surgeries, where each leaf of each clasper has total linking number
zero with the link.

iii. $L_{0}$ and $L_{1}$ are Strongly S-equivalent.

iv. For some choice of Seifert Surfaces $\Sigma _{0}$ and $\Sigma
_{1}$ and bases of $H_{1}(\Sigma _{i})$, $L_{0}$ and $L_{1}$ have
identical ordered Seifert Matrices.
\end{main}

Throughout the paper, a\emph{\ link} $L$ with $m$-components is a subset of $%
S^{3}$, or of $\mathbb{R}^{3}$, that consists of $m$ disjoint,
piecewise linear, simple closed curves (a link with one component is
a \emph{knot}). \ Unless otherwise noted, all links will be oriented
and ordered. \ Seifert surfaces for links will be required to be
connected, and $g$ will denote the minimal genus among possible
Seifert surfaces.

\section{The Doubled-Delta Move}

\subsection{The Delta Move}

\begin{figure}[h]
 \centering
 \includegraphics{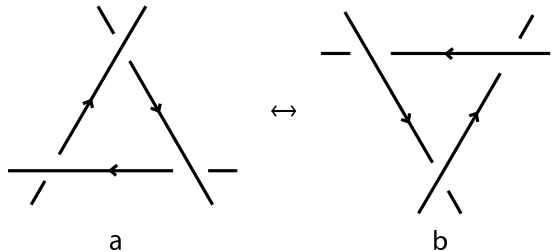}
 \caption{The Delta Move}
 \label{fig delta}
\end{figure}

Before defining the doubled-delta move, we look first at the simpler
(single) delta move. \ The \emph{delta move} shown in Figure \ref{fig delta}%
\ is a particular move on knot or link diagrams. \ Given a link containing
the tangle in Figure \ref{fig delta}a, replace the tangle with that of
Figure \ref{fig delta}b in such a way that respects the free ends.

It is clear that the delta move preserves the pairwise linking
numbers of the components, since for each pair of strands, the move
changes neither the crossing nor the orientation. \ Though less
obvious, the converse is also true, as proved by Murakami and Nakanishi in \cite{Murakami}%
: \ two links have the same sets of pairwise linking numbers if and
only if they are equivalent under delta moves. \ Since a knot is a
one-component link, it follows that any two knots are equivalent
under delta moves.

\subsection{The Doubled-Delta Move}

The \emph{doubled-delta move}\ is\emph{\ }similar to the single
delta move, with each of the three strands being replaced by a pair
of oppositely oriented strands. \ The link components that
comprise the six affected strands are irrelevant. \ The
doubled-delta move is more restrictive than the (single) delta move,
but it still preserves the pairwise linking numbers.

\begin{figure}[h]
 \centering
 \includegraphics{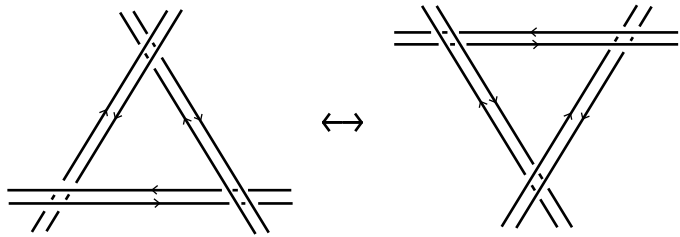}
 \caption{The Doubled-Delta Move}
 \label{fig ddmove1}
\end{figure}

\subsection{Borromean Surgery\label{section borromeo surgery}}

\begin{figure}[h]
 \centering
 \includegraphics{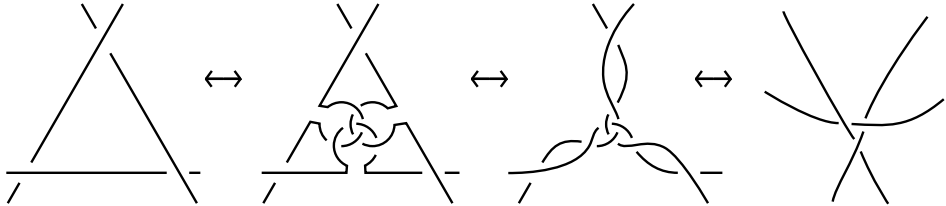}
 \caption{The Delta-Move is Borromean Surgery}
 \label{fig delta=borr}
\end{figure}

The delta move and doubled-delta move partition knots and links into
equivalence classes. \ Two links are said to be in the same class if they
are related by a sequence of such moves. \ However, unlike many of the basic
moves that change a knot or link diagram, the delta and doubled-delta move
are closely related to other operations used by topologists in a variety of
contexts. \ In Figure \ref{fig delta=borr}, we see that the effect of the
delta move is the same as \textquotedblleft adding in\textquotedblright\
Borromean rings.

Another name for the operation depicted in the figure is Borromean \
(or Borromeo) surgery. \ We will see in Section \ref{section
clasper} how Borromean surgery relates to claspers, and thus to
grope cobordisms, finite-type invariants, and other ideas of
low-dimensional topology.

\section{Claspers\label{section clasper}}

Under what conditions are two links related by a sequence of doubled-delta
moves? \ This question is interesting to many topologists in the context of
clasper surgery, a special case of which is Borromean surgery. \ Claspers were
first defined by K. Habiro \cite{Habiro}, where they arose in the context of
finite-type invariants. \ Habiro demonstrated how the theory of claspers
provides an alternative calculus under which one can study finite-type
invariants of knots and 3-manifolds. Claspers also were implicit in the work
of Goussarov \cite{Goussarov1}\cite{Goussarov2}. \ Today, they are studied
across various fields of low-dimensional topology. \ In addition to
applications of finite-type invariants, P.\ Teichner and J. Conant examine
claspers' relationships to grope cobordisms \cite{conant}, \ and S.
Garoufalidis uses clasper surgery \cite{stavros2} to better understand the
Kontsevich integral and concordance classes of knots. \

\subsection{Definition}

A \emph{clasper }is a compact surface constructed from the following three
types of pieces:

\begin{itemize}
\item \emph{edges}, or bands that connect the other two types of pieces

\item \emph{nodes}, or disks with three incident edges

\item \emph{leaves}, or annuli with one incident edge. \
\end{itemize}

The annuli that comprise the leaves may be twisted with any number of full
twists. \ We call this number the \emph{framing} of the leaf. \ Figure \ref%
{fig clasperlink}a shows a clasper with zero-framed leaves.

\begin{figure}[h]
 \centering
 \includegraphics{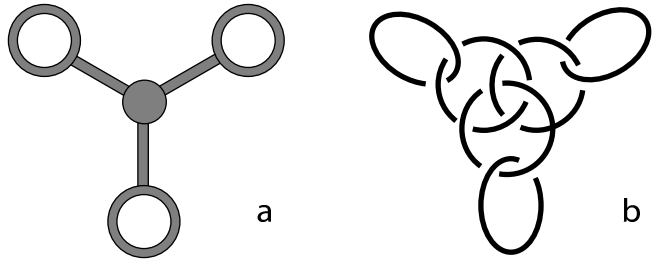}
 \caption{A $Y$-clasper, with its associated link.}
 \label{fig clasperlink}
\end{figure}

\subsection{Clasper Surgery}

Assume that a clasper $C$ is embedded in a 3-manifold $M$. \ To do surgery
on the clasper means to remove a handlebody neighborhood of the clasper from
a $3$-manifold $M^{3}$, and glue it back in a prescribed way according to
the clasper and its framing. \ In particular, we associate a link $L_{C}$ to
the clasper $C$ as in Figure \ref{fig clasperlink}b using the following
substitutions: \ each node of the clasper is replaced by a copy of the
positive, zero-framed Borromean rings, and each edge is replaced by a
positive Hopf link. \ Leaves of the clasper do not contribute additional
link components, but the framing of a leaf does determine the framing of the
Hopf link component corresponding to that leaf's adjacent edge. \ Clasper
surgery, then, is really integer surgery on the associated framed link.

An important property of clasper surgery is that it preserves the homology
of the affected manifold. \ When each leaf of the clasper \textquotedblleft
clasps\textquotedblright\ a knot or link and at least one leaf bounds a disk
in $M^{3}-clasper$, the result of the surgery is a new knot or link in the
same $3$-manifold.

Surgery on the simplest clasper--the strut, with a single edge and
two leaves--can accomplish a single crossing change in the original
knot. \ Since the crossing change is an unknotting relation,
all knots are related by strut-clasper surgery. \ The simplest interesting clasper, then, is the $%
Y $-clasper. \ One may also argue that the $Y$-clasper is the most
interesting clasper, since any larger clasper may be realized as several $Y$%
-claspers by expanding edges of the larger clasper into Hopf-linked pairs of
leaves \cite{Habiro}. \

\subsection{Null Claspers}

In \cite{stavros3}, S. Garoufalidis and L. Rozansky discuss \emph{null
claspers,} those whose leaves are null-homologous links in $M-K$ for a
3-manifold $M$ and knot $K$. \ Such leaves are sent to zero under the map $%
\pi _{1}(M-K)\rightarrow H_{1}(M-K)\cong \mathbb{Z}$; \ in other
words, each leaf has algebraic linking number zero with the knot. \
Garoufalidis and Rozansky explain that null claspers have been used to demonstrate a
rational version of the Kontsevich integral and can be used to
define a notion of finite-type invariants. \ In Lemma 1.3 of
\cite{stavros3}, they show that surgery on null claspers preserves
not only the homology of a knot complement, but also the Alexander
module and Blanchfield linking form. \ Furthermore, null clasper
surgery describes a move on the set of knots in integral homology
spheres that directly corresponds to the doubled-delta move. \

Extending Garoufalidis's and Rozansky's definition to links, for the pair $(M,L)$, 
where $M$ is a 3-manifold and $L$ is a link, a null
clasper would be one whose leaves have linking number zero with each
component. \ Unfortunately, the doubled-delta move acts on the strands
of a link independently of the link components. \ Therefore, in considering the
relationship between claspers and the doubled-delta move, our interest is
with a slightly larger class of claspers: $\ $those in which each leaf
clasps several strands of the link in such a way that the \emph{total}
linking number with all link components is zero. \

We claim that when the leaves of a zero-framed $Y$-clasper each have \emph{%
total} linking number zero with the link, the $Y$-clasper surgery has the
same effect on the link as a finite sequence of doubled-delta moves, as well
as that of Borromean surgery. \ This is one of the implications of our main
theorem, and is illustrated in Figure \ref{fig ddelta=clasper}.

\section{S-equivalence and the Doubled-Delta Move}

\subsection{The Seifert Matrix}

Many of the basic algebraic invariants of knot theory, including\ the
Alexander module, the Conway polynomial, and the signature, depend only on
the Seifert matrix of a knot or link, which is readily computable but not
uniquely defined. \

For every $m$-component link $L$ in $S^{3}$, there is a Seifert surface $%
\Sigma $ associated to the link, where $\Sigma $ is a connected, oriented,
embedded surface with the components of $L$ as its boundary. \ Given a basis
$\{b_{i}\}$ of $H_{1}(\Sigma )$, we can associate a \emph{Seifert matrix} $M$
to the link $L$, where the entries of $M$ are defined from the linking
number of two basis elements. \ In particular, $M_{i,j}=lk(b_{i},b_{j}^{+})$%
, where $b_{j}^{+}$ is a pushoff of $b_{j}$ in the positive normal
direction. \

\subsection{Classical S-equivalence}

S-equivalence is a notion that has been widely considered for both knots and
links \cite{gordon} \cite{Kawauchi} \cite{lickorish}. \ Two square integral
matrices $M$ and $N$ are said to be \emph{S-equivalent} if $M$ can be
transformed into $N$ by a finite sequence of integral congruences (that is, $%
M=A^{t}NA$ for some integral matrix $A$ with $\det (A)=\pm 1$) and row or
column enlargements/reductions of the form

\begin{equation*}
N=%
\begin{pmatrix}
&  &  & | & | \\
& M &  & \overrightarrow{y}^{t} & 0 \\
&  &  & | & | \\
\text{---} & \overrightarrow{x} & \text{---} & z & 0 \\
\text{---} & 0 & \text{---} & 1 & 0%
\end{pmatrix}%
\ \ \ or\ \ N=%
\begin{pmatrix}
&  &  & | & | \\
& M &  & \overrightarrow{y}^{t} & 0 \\
&  &  & | & | \\
\text{---} & \overrightarrow{x} & \text{---} & z & 1 \\
\text{---} & 0 & \text{---} & 0 & 0%
\end{pmatrix}%
.
\end{equation*}

A more precise statement of these definitions is given in section \ref{section sseq}.  Up to S-equivalence, Seifert matrices are well-defined for knots and links.

\subsection{Relationship to the Doubled-Delta Move}

Despite being a purely algebraic relation, S-equivalence has geometric
implications for knots. \ In particular, S. Naik and T. Stanford prove that
two knots are S-equivalent if and only if they are equivalent by a sequence
of doubled-delta moves \cite{ns}. \

One might assume that a similar statement should be true of links. \ Is the
S-equivalence of two links enough to guarantee that they are equivalent
under a sequence of doubled-delta moves? \ This question was posed by
Stavros Garoufalidis in the context of clasper surgery, and given that all
the definitions leading up to S-equivalence are the same for links as they
are for knots, it seems that the analogous proof for links should follow in
a straightforward fashion from the proof of Naik-Stanford. \ The question
seems plausible, but in fact it is false. \ We prove by counterexample in
Proposition \ref{prop SSnecessary} that two S-equivalent links are not
necessarily related by a sequence of doubled-delta moves.

\begin{proposition}
\label{prop SSnecessary}S-equivalence is not a sufficient condition for two
links to be related by a sequence of doubled-delta moves.
\end{proposition}

\begin{proof}
The two links $L_{0}$ and $L_{1}$ depicted in figure \ref{fig
SnotSSlinkex} are S-equivalent but not related by a sequence of
doubled delta moves.
\begin{figure}[h]
 \centering
 \includegraphics{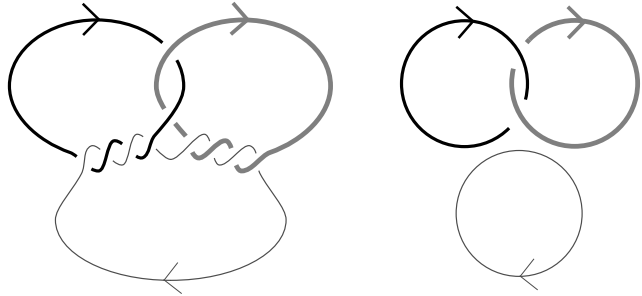}
 \caption{Two 3-component Links}
 \label{fig SnotSSlinkex}
\end{figure}
Note that the pairwise linking numbers of the first link are
$\left\{ -1,2,2\right\}$ while for the second link they are $\left\{
1,0,0\right\}$. Since the doubled-delta move preserves pairwise
linking numbers, these two links cannot be related by doubled-delta
moves.

The fact that the two links are S-equivalent can be seen by choosing Seifert
Surfaces as in figure \ref{fig surface1ab},
\begin{figure}[h]
 \centering
 \includegraphics{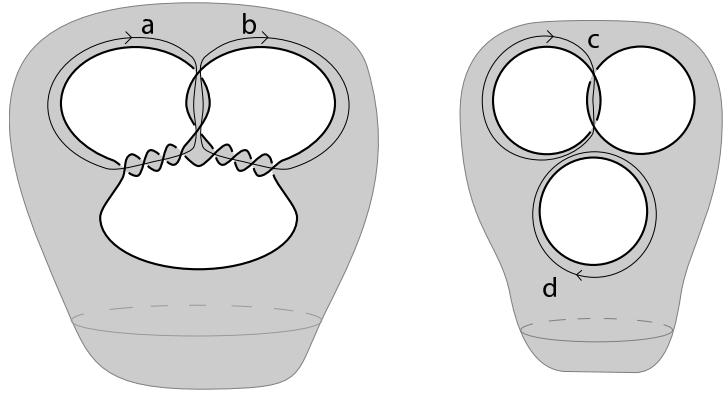}
 \caption{Two Seifert Surfaces}
 \label{fig surface1ab}
\end{figure}
with bases $\{a,b\}$ and $\{c,d\}$. \ Each of the surfaces can be
viewed as a punctured 2-sphere with a Y-shaped band glued along the
boundary, where each of the three strips of the Y is twisted
according to the desired linking numbers. \ The Seifert matrices are
then $M_{0}=\left(
\begin{smallmatrix}
-1 & -1 \\
-1 & -1%
\end{smallmatrix}%
\right) $ and $M_{1}=\left(
\begin{smallmatrix}
-1 & 0 \\
0 & 0%
\end{smallmatrix}%
\right) $, respectively. \ The matrix $A=\left(
\begin{smallmatrix}
1 & 1 \\
0 & 1%
\end{smallmatrix}%
\right) $ satisfies the condition $A^{t}M_{1}A=M_{0}$, demonstrating that
the two links are in fact S-equivalent.
\end{proof}

Finding the sufficient conditions to extend the Naik-Stanford theorem will
require some new definitions, treated in section \ref{section sseq}. \

\section{Strong S-equivalence\label{section sseq}}

As shown in Proposition \ref{prop SSnecessary}, the classical definition of
S-equivalence is inadequate for the proof of our main theorem, and actually
it seems inappropriate for links in many contexts. \ Treating a link as a
disjoint knot, S-equivalence regards the entire link as a whole, without
reference to the individual components. \ Strong S-equivalence, defined
below, better respects the boundary components of a Seifert surface, and
hence the components of the link. \

Let $L=\{L_{1},L_{2},...,L_{m}\}$ be an oriented, ordered link in $S^{3}$. \
Let $\Sigma $ be a Seifert surface for $L$ with $m$ boundary components, and
let $g$ denote the genus of $\Sigma $. \ That is, $\Sigma \subseteq S^{3}$
is an oriented surface with $\partial \Sigma =L$. \ Construct an \emph{%
ordered basis} $\beta =\{\ell _{1},\ell _{2},...\ell _{m-1},\beta _{1},\beta
_{2,}...\beta _{2g}\}$ for $H_{1}(\Sigma )$, where $\ell _{i}$ is
represented by the $i$th component $L_{i}$ of $L$. \ Define the Seifert
pairing $\sigma (a,b)=lk(a,b^{+})$, where $b^{+}$ is a pushoff of $b$ in the
positive normal direction. \ We introduce the following term:

\begin{definition}
A matrix $M$ representing $\sigma $ with respect to an ordered basis $\beta $
of the form above is called an \emph{Ordered Seifert Matrix} for the
oriented ordered link $L$. \
\end{definition}

\begin{remark}
\label{remark lambda block}$M$ has the form of a block matrix $\left(
\begin{smallmatrix}
\lambda & A \\
B & C%
\end{smallmatrix}%
\right) $, where $\lambda $ is an $(m-1)\times (m-1)$ block and $C$ is a $%
2g\times 2g$ block. \ The block $\lambda $ is completely determined by the
pairwise linking numbers of $L$, and has the following properties:

$i.$ \ $\lambda _{i,i}=-\sum\limits_{i\neq j}\lambda _{i,j}-lk(L_{i},L_{m})$
for $1\leq i,j\leq m-1$, and

$ii.$ \ $\lambda _{i,j}=\lambda _{j,i}=lk(L_{i},L_{j})$ if $i\neq j$ and $%
1\leq i,j\leq m-1$.
\end{remark}

The first of these properties is somewhat less obvious. \ If we let $\ell
_{m}$ be the homology class represented by $L_{m}$, then the class $%
\sum\limits_{j=1}^{m}\ell _{j}=0\in H_{1}(\Sigma )$, since it is represented
by the boundary of $\Sigma $. \ Then
\begin{eqnarray*}
0 &=&\sigma \left( \ell _{i},\sum\limits_{j=1}^{m}\ell _{j}\right) \\
&=&\sum\limits_{j=1}^{m}\sigma (\ell _{i},\ell _{j}) \\
&=&\sum\limits_{i\neq j}\sigma (\ell _{i},\ell _{j})+\sigma (\ell _{i},\ell
_{m})+\sigma (\ell _{i},\ell _{i})\text{ \ \ \ for }1\leq i,j\leq m-1 \\
&=&\sum\limits_{i\neq j}\lambda _{i,j}+lk(L_{i},L_{m})+\lambda _{i,i}\text{
\ \ \ for }1\leq i,j\leq m-1\text{.}
\end{eqnarray*}
\

The following two definitions can be found in \cite{Kawauchi}, as they are
essential to the classical definition of S-equivalence. \

\begin{definition}
We say two integral square matrices $V$ and $W$ are \emph{congruent} if $%
V=P^{t}WP$ for some integral matrix $P$ with $\det (P)=\pm 1$. \
\end{definition}

\begin{definition}
For integral square matrices $V$ and $W$, we say that $W$ is an \emph{%
enlargement} of $V$, or $V$ is a \emph{reduction} of $W$ if

\begin{equation*}
W=%
\begin{pmatrix}
&  &  & | & | \\
& V &  & \overrightarrow{y}^{t} & 0 \\
&  &  & | & | \\
\text{---} & \overrightarrow{x} & \text{---} & z & 1 \\
\text{---} & 0 & \text{---} & 0 & 0%
\end{pmatrix}%
\text{ \ \ \ or \ \ }W=%
\begin{pmatrix}
&  &  & | & | \\
& V &  & \overrightarrow{y}^{t} & 0 \\
&  &  & | & | \\
\text{---} & \overrightarrow{x} & \text{---} & z & 0 \\
\text{---} & 0 & \text{---} & 1 & 0%
\end{pmatrix}%
.
\end{equation*}
\end{definition}

With a slight modification of these two definitions, we introduce:

\begin{definition}
Two matrices $V$ and $W$ are \emph{Strongly S-equivalent} if $V$ is
equivalent to $W$ under a finite sequence of

\begin{itemize}
\item Congruences $(\cong )$ that fix the upper-left $(m-1)\times (m-1)$
block of the matrix. \ That is, there exists an integral matrix $A=\left(
\begin{smallmatrix}
I & \ast \\
0 & \ast%
\end{smallmatrix}%
\right) $, so that $A^{t}VA=W$ and $\det (A)=\pm 1$, where $I$ is the $%
(m-1)\times (m-1)$ identity matrix. \

\item Enlargements $(\nearrow )$ and reductions $(\searrow )$ where the
first $(m-1)$ elements of the vectors $\overrightarrow{x}$ and $%
\overrightarrow{y}$ are equal, and where the reduced matrix $V$ is $n\times
n $ for $n\geq m-1$. \
\end{itemize}
\end{definition}

Note that reductions should not be allowed to reduce the size of the matrix
smaller than $(m-1)\times (m-1)$, for in the case of Seifert matrices, this
would effectively eliminate link components. \

Also note that the first $(m-1)$ elements of the vectors $\overrightarrow{x}$
and $\overrightarrow{y}$ are equal because these entries should be zero in
the intersection form $W-W^{t}$ of the enlarged matrix $W$, since boundary
components will not intersect any other basis elements. \

\begin{definition}
We say that two links $L$ and $L^{\prime }$ are \emph{Strongly S-equivalent}
if, for some choice of Seifert surfaces and ordered bases, $L$ and $%
L^{\prime }$ have ordered Seifert matrices $M$ and $M^{\prime }$ that are
Strongly S-equivalent .
\end{definition}

One might object that Strong S-equivalence imposes the \textquotedblleft
restriction\textquotedblright\ that $V$ and $W$ agree on their upper-left $%
(m-1)\times (m-1)$ blocks. \ However, since homeomorphisms of $\Sigma $ from
the pure mapping class group preserve the boundary components pointwise, any
change of basis of $M$ must fix the $\ell _{1},\ell _{2},...\ell _{m-1}$
basis elements. \ With Strong S-equivalence, we're not \textquotedblleft
restricting\textquotedblright\ our definition so much as \textit{respecting}
the boundary components of a link. \

\subsection{Disk-Band Form of the Seifert Surface}

It will often be useful to look at the Seifert surface in disk-band
form. \ Any Seifert surface for a link is homeomorphic to a disk
with $2g+m-1$ bands attached; more importantly, we can always
isotope our surface to one in disk-band form as in Figure \ref{fig
diskband stringlink w/basis}, with $2g$ bands interlaced in pairs
and $m-1$ non-interlaced bands.
\begin{figure}[h]
 \centering
 \includegraphics{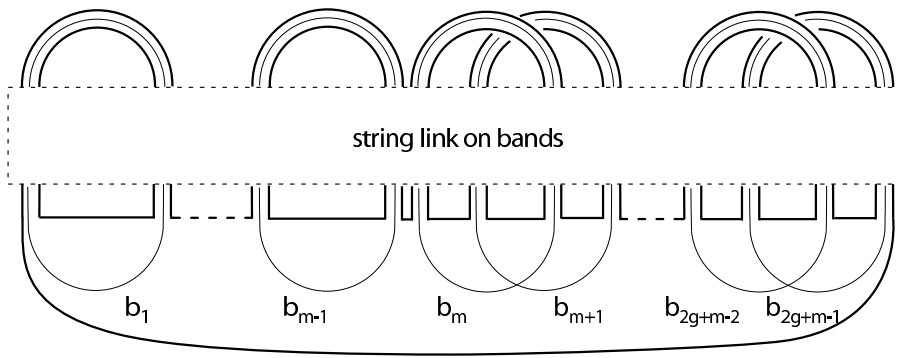}
 \caption{Disk-Band Form for a Link's Seifert Surface}
 \label{fig diskband stringlink w/basis}
\end{figure}
Within the dotted box, the bands of our surface form a \emph{string
link}, i.e. they may be knotted, twisted, or intertwined as long as
the bands entering the top of the box match up with those that leave
the bottom of the box. \ With the
surface in this form, there is an natural choice for a basis $\{b_{i}\}$ of $%
H_{1}(\Sigma )$, with one basis element running through each of the bands. \
The first $m-1$ non-interlacing bands yield basis elements that are pushoffs
of the respective boundary components and have no intersection with any
other basis elements. \ The $2g$ interlaced bands yield basis elements with
the properties of a symplectic basis, i.e. one whose intersection form $%
V_{i,j}=\langle b_{i},b_{j}\rangle $ is the block sum $\bigoplus%
\limits_{j=1}^{g}\left(
\begin{smallmatrix}
0 & -1 \\
1 & 0%
\end{smallmatrix}%
\right) $. \ We coin the term \emph{semi-symplectic}\ to refer to this
natural choice of basis. \ \ A semi-symplectic basis is an ordered basis,
but the converse is not necessarily true.

\subsection{Properties of Strong S-equivalence}

\begin{proposition}
Pairwise linking number is an invariant of Strong S-equivalence.
\end{proposition}

\begin{proof}
If two links $L$ and $L^{\prime }$ are Strongly S-equivalent then,
regardless of the choice of ordered Seifert matrices $M$ and $M^{\prime }$,
the upper left $(m-1)\times (m-1)$ blocks of $M$ and $M^{\prime }$ will
necessarily agree. \ That is, for $i\neq j$, $1\leq i,j<m$,
\begin{equation*}
lk(L_{i},L_{j})=M_{i,j}=M_{i,j}^{\prime }=lk(L_{i}^{\prime },L_{j}^{\prime })%
\text{.}
\end{equation*}%
Furthermore, by property $i.$ of Remark \ref{remark lambda block},%
\begin{equation*}
lk(L_{i},L_{m})=-\sum\limits_{j=1}^{m-1}M_{i,j}=-\sum%
\limits_{j=1}^{m-1}M_{i,j}^{\prime }=lk(L_{i}^{\prime },L_{m}^{\prime }).
\end{equation*}%
\ Thus all the pairwise linking numbers for $L$ agree with those for $%
L^{\prime }$. \
\end{proof}

\begin{proposition}
\label{prop M1=M2}Any two ordered Seifert matrices for an oriented ordered
link\ $L$ are Strongly S-equivalent. \
\end{proposition}

\begin{proof}
Let $M_{1}$ and $M_{2}$ be two ordered Seifert matrices for $L$ with respect
to Seifert surfaces $\Sigma _{1}$ and $\Sigma _{2}$ and bases $\beta _{1}$
and $\beta _{2}$. \

By Lemma 5.2.4 of \cite{Kawauchi}, any two connected Seifert surfaces $%
\Sigma _{1}$ and $\Sigma _{2}$ of a link $L$ are ambient isotopic after
modifying them by a finite sequence of 1-handle enlargements. \ In other
words, there are two ambient isotopic surfaces $\widetilde{\Sigma _{1}}$ and
$\widetilde{\Sigma _{2}}$ such that $\widetilde{\Sigma _{i}}$ is the result
of several 1-handle enlargements of $\Sigma _{i}$. \ Without loss of
generality, we can consider just the following case: \ suppose some surface $%
\widehat{\Sigma _{2}}$ is ambient isotopic to $\widehat{\Sigma _{1}}$, a
single 1-handle enlargement of $\Sigma _{1}$. \ Let $a_{2}$ be a meridian\
(technically, the belt sphere) of the 1-handle and choose a closed curve $%
a_{1}$ that intersects the meridinal disk (technically, the co-core) exactly
once. \ Then one of

\begin{equation*}
\widehat{M_{1}}=%
\begin{pmatrix}
&  &  & | & | \\
& M_{1} &  & \overrightarrow{y}^{t} & 0 \\
&  &  & | & | \\
\text{---} & \overrightarrow{x} & \text{---} & z & 1 \\
\text{---} & 0 & \text{---} & 0 & 0%
\end{pmatrix}%
\text{ \ \ or \ \ }\overline{M_{1}}=%
\begin{pmatrix}
&  &  & | & | \\
& M_{1} &  & \overrightarrow{y}^{t} & 0 \\
&  &  & | & | \\
\text{---} & \overrightarrow{x} & \text{---} & z & 0 \\
\text{---} & 0 & \text{---} & 1 & 0%
\end{pmatrix}%
\end{equation*}%
is a Seifert matrix for $\widehat{\Sigma _{1}}$ with respect to the basis $%
\beta _{1}\cup \{a_{1},a_{2}\}$ for $H_{1}(\widehat{\Sigma _{1}})$. \ The
matrix on the left, $\widehat{M_{1}}$, corresponds to a 1-handle attached so
that its core lies on the negative side of the surface, while $\overline{%
M_{1}}$ corresponds to a 1-handle with its core on the positive side of the
surface. \ We will treat only the $\widehat{M_{1}}$ case and note that the $%
\overline{M_{1}}$ case follows similarly. \ The zeros in the last column of
$\widehat{M_{1}}$ come from the fact that a positive pushoff of the
meridian, $a_{2}$, will not link any of the basis elements except $a_{1}$. \
The last row of $\widehat{M_{1}}$ is all zero because a negative pushoff of
the meridian $a_{2}$ will not link any of the basis elements. \ Since for
any Seifert matrix $N$, $N-N^{t}$ is an intersection form, the first $m-1$
entries of the vectors $\overrightarrow{x}$ and $\overrightarrow{y}$ must be
equal. \ Otherwise $a_{1}$ would be intersecting a boundary component, which
is impossible. \ The rest of the entries of $\overrightarrow{x}$ and $%
\overrightarrow{y}$ are freely determined by the particular embedding of the
new 1-handle and the choice of curve $a_{1}$, as always with $%
M_{i,j}=lk(\beta _{i},\beta _{j}^{+})$. \ The entry $z$ is, of course, $%
lk(a_{1},a_{1}^{+})$. \ Note that with more restrictions on the choice of $%
a_{1}$, the allowable entries of $\overrightarrow{x}$, $\overrightarrow{y}$,
and $z$ could be given more structure. \

By the definitions of enlargement and reduction, $\widehat{M_{1}}$ and $M_{1}$
are Strongly S-equivalent. \ This process can be iterated for any finite
number of 1-handle enlargements so that the final matrix $\widetilde{M_{1}}$ is Strongly
S-equivalent to $M_{1}$, and similarly, $M_{2}$ can be shown to be Strongly
S-equivalent to $\widetilde{M_{2}}$, the similarly constructed Seifert
matrix for $\widetilde{\Sigma _{2}}$. \

Since $\widetilde{\Sigma _{1}}$ is ambient isotopic to $\widetilde{\Sigma
_{2}}$, $\widetilde{M_{1}}$ and $\widetilde{M_{2}}$ differ only by a choice
of bases. \ Furthermore, since the upper left blocks of both $\widetilde{%
M_{1}}$ and $\widetilde{M_{2}}$ are determined by the link $L$ and are equal
(note that $\widetilde{M_{i}}$ has the same upper-left block as $M_{i}$), $%
\widetilde{M_{1}}$ and $\widetilde{M_{2}}$ are related by a change of basis
that preserves the first $m-1$ basis elements, and they are thus Strongly
S-equivalent. \ By transitivity, $M_{1}$ is Strongly S-equivalent
to $M_{2}$. \
\end{proof}

Like its classical S-equivalence analogue, the converse of Proposition \ref%
{prop M1=M2} is not true. \ If $N$ is a Seifert matrix for the link $L$ with
respect to Seifert surface $\Sigma $, and $M$ is Strongly S-equivalent to $N$%
, $M$ is not necessarily a Seifert matrix for $L$. \ The difficulty arises
when $M$ is a reduction of $N$, since it may not be possible to find a
corresponding 1-handle reduction of the Seifert surface $\Sigma $. \ The
converse \emph{is} true, however, in the special cases of congruences and
enlargements. \ Matrix congruence corresponds to change of basis for widecheck$%
H_{1}(\Sigma )$. \ In order to explicitly understand how any given matrix
enlargement corresponds to a 1-handle enlargement of a Seifert surface, we
prove the following proposition: \

\begin{proposition}
\label{prop enlargement}If $N$ is an ordered Seifert matrix for the link $L$
with respect to Seifert surface $\Sigma $, and $M$ is an enlargement of $N$,
then $M$ is an ordered Seifert matrix for $L$ with respect to a 1-handle
enlargement $\widehat{\Sigma }$ of $\Sigma .$
\end{proposition}

\begin{proof}
Let $\Sigma $ be a Seifert surface for $L$ with basis $\beta $ of 
$H_{1}(\Sigma )$ that induces the Seifert matrix $N$.  We will construct 
the 1-handle enlargement $\widehat{\Sigma }$ of $\Sigma $
that has $M$ as its Seifert matrix, where either

\begin{equation*}
M=%
\begin{pmatrix}
&  &  & | & | \\
& N &  & \overrightarrow{y}^{t} & 0 \\
&  &  & | & | \\
\text{---} & \overrightarrow{x} & \text{---} & z & 1 \\
\text{---} & 0 & \text{---} & 0 & 0%
\end{pmatrix}%
\text{ \ or \ }M=%
\begin{pmatrix}
&  &  & | & | \\
& N &  & \overrightarrow{y}^{t} & 0 \\
&  &  & | & | \\
\text{---} & \overrightarrow{x} & \text{---} & z & 0 \\
\text{---} & 0 & \text{---} & 1 & 0%
\end{pmatrix}%
.
\end{equation*}

As in Proposition \ref{prop M1=M2},
the matrix on the left corresponds to attaching the 1-handle so that its
core lies on the negative side of the surface, while the matrix on the right
corresponds to a 1-handle with its core on the positive side of the surface.
\

Let $\{b_{1},\ldots ,b_{2g+m-1}\}\subset \Sigma $ be a set of representative
curves for the basis $\beta $. \ Find two small disks in $\Sigma $ that are
disjoint from the curves $\{b_{1},\ldots ,b_{2g+m-1}\}$. \ These disks will
be the attaching region for the 1-handle. \ Designate two points, $p$ and $q$%
, one point on the boundary of each disk. \  The enlarged matrix $M$ determines 
the Seifert form--that is, a combination of linking and intersection numbers--of 
two new basis elements $a_{1}$ and $a_{2}$ with
curves $\{b_{1},\ldots ,b_{2g+m-1}\}$. \ The 1-handle we construct will have
$a_{2}$ as a meridian and $a_{1}$ running parallel to the core. \ First let
us construct $a_{1}$, breaking it into two parts, $\gamma $ and $\delta $, as
in figure \ref{fig 1handle}. In order to be consistent with the Seifert
matrix $M$, we need $a_{1}$ to to intersect and link the curves $%
\{b_{1},\ldots ,b_{2g+m-1}\}$ according to the entries of the vectors $%
\overrightarrow{x}$ and $\overrightarrow{y}$. \ We can extract the
intersection information from the intersection form $M-M^{t}$, or in
particular, from the vector entries $\overrightarrow{x_{i}}-\overrightarrow{%
y_{i}}$ \cite{Rolfsen}. \ Running through the surface $\Sigma $, $\gamma $
will be constructed to take care of any intersections, while the handle itself, and
hence $\delta $, will be free to link the basis elements. \ Together, $%
\gamma \cup \delta =a_{1}$ will then satisfy each of the entries $%
lk(a_{1},b_{i}^{+})=x_{i}$ and $lk(b_{i},a_{1}^{+})=y_{i}$.

\begin{figure}[h]
 \centering
 \includegraphics{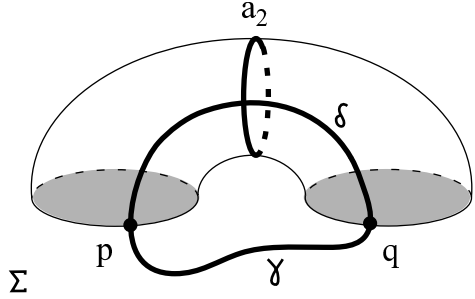}
 \caption{1-handle}
 \label{fig 1handle}
\end{figure}

To find $\gamma \subset \Sigma $, the half of $a_{1}$ that will travel from $%
p$ to $q$ along the surface $\Sigma $, let us look at the intersection form $%
M-M^{t}$. \ The individual entries of $\overrightarrow{x}-\overrightarrow{y}$
determine how $\gamma $ should intersect the curves $\{b_{1},\ldots
,b_{2g+m-1}\}$; since $\delta \subset a_{1}$ does not intersect $\beta $ at
all, $\langle \gamma ,b_{i}\rangle =\langle a_{1},b_{i}\rangle =%
\overrightarrow{x_{i}}-\overrightarrow{y_{i}}$. \ The first $m-1$ entries of
$\overrightarrow{x}-\overrightarrow{y}$ are zero, which is consistent with
the fact that $\gamma $ cannot cross a boundary component. \ For the last $%
2g $ curves $b_{i}$, choose $\gamma $ so that $\langle \gamma ,b_{i}\rangle =%
\overrightarrow{x_{i}}-\overrightarrow{y_{i}}$. \ This is possible. \ As
described previously, the basis $\beta $ is in correspondence with a
semi-symplectic basis. \ For ease of construction, choose a set of curves $%
\{c_{i}\}$ in $\Sigma $ to be representatives of this semi-symplectic basis\
with the same (semi-symplectic) algebraic intersection properties. \ If, for any set of
integers $\{k_{i}\}$, we can construct a curve $\gamma $ that has algebraic
intersection number $k_{i}$ with each $c_{i}$, then the same can be done for
the set of curves $\left\{ b_{i}\right\} $ via this correspondence.\ Start
with a curve $\gamma _{1}$ from $p$ to $q$ that does not intersect any of
the curves $c_{i}$. \ For each desired $(\pm )$ intersection with a
particular curve $c_{2k}$, take $\gamma _{2}=\gamma _{1}\pm c_{2k-1}.$ \
(For each $(\pm )$ intersection with $c_{2k-1}$, take $\gamma _{2}=\gamma
_{1}\mp c_{2k}$.) Continue this until all intersections are achieved, and
call the final curve $\gamma $. \ Figure \ref{fig handleconstruct}, below,
demonstrates how a path $\gamma $ from $p$ to $q$ can be chosen to intersect
$c_{1}$ exactly one time and intersect $c_{4}$ exactly $-2$ times.

\begin{figure}[h]
 \centering
 \includegraphics{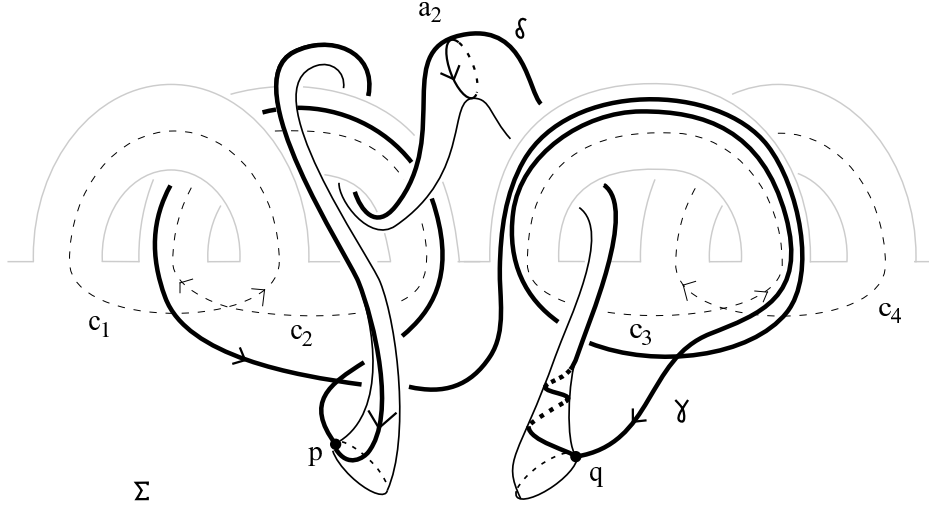}
 \caption{Construction of 1-handle}
 \label{fig handleconstruct}
\end{figure}

After the intersection properties are established between $\gamma $ and the $%
b_{i}$, we shift our focus to the linking properties. \ As with the
intersections, the linking numbers of $\gamma $ with the $b_{i}$ are in
correspondence with those of $\gamma $ with the $c_{i}$, so we will focus
our construction on the semi-symplectic basis $\{c_{i}\}$. \ The core of the
1-handle can be chosen in the complement of $\Sigma $ so that it links each
of the curves $c_{i}$ the desired number of times. \ This is easy to see in
the disk-band representation of $\Sigma $. \ The core can then be fattened
up to a solid handle, the surface of which is the 1-handle enlargement of $%
\Sigma .$ \ We choose a curve $\delta $ on the surface of the handle
parallel to the core, and let $a_{1}=\gamma +\delta $ be one of our new
basis elements. \ Figure \ref{fig handleconstruct} demonstrates how the core
can be chosen to link $c_{2}$ one time and link $c_{3}$ negative one time. \

Lastly, we must adjust the path $\delta $ along the 1-handle so that $%
lk(a_{1},a_{1}^{+})=z$. \ The remedy is simple: \ each repeated positive or
negative full twist to the 1-handle will increase or decrease $%
lk(a_{1},a_{1}^{+})$ by one without affecting any of the other basis
elements. \ Figure \ref{fig handleconstruct} illustrates $\delta $ chosen so
that $z=-2.$ \

We have explicitly constructed a new Seifert surface $\widehat{\Sigma }$ for
$L$ and a new basis $\{\beta ,a_{1},a_{2}\}$ for $H_{1}(\widehat{\Sigma })$ 
such that $M$ is an
ordered Seifert matrix for $L$ with respect to them. \
\end{proof}

\section{Proof of the Main Theorem}

\begin{theorem}
\label{mainthm}Consider two oriented, ordered $m$-component links $L_{0}$
and $L_{1}$. \ The following four statements are equivalent:

i. $\ L_{1}$ can be obtained from $L_{0}$ through a sequence of
doubled-delta moves.

ii. $\ L_{0}$ and $L_{1}$ are related by a sequence of $Y$-clasper
surgeries, where each clasper has total linking number zero with the link. \

iii. \ $L_{0}$ and $L_{1}$ are Strongly S-equivalent.

iv. \ For some choice of Seifert Surfaces $\Sigma _{0}$ and $\Sigma _{1}$
and bases of $H_{1}(\Sigma _{i}),$ $L_{0}$ and $L_{1}$ have the same ordered
Seifert Matrix.
\end{theorem}

The proofs of each implication $i.\Longrightarrow ii.\Longrightarrow
iii.\Longrightarrow iv.\Longrightarrow i.$ are treated individually below:

\subsection{Proof of $i.\Longrightarrow ii.$}

\begin{proof}
Doubled-delta moves correspond to \textquotedblleft Borromean
surgery,\textquotedblright\ which is exactly the effect of
$Y$-clasper surgery, where each leaf clasps pairs of oppositely
oriented strands. \ This is depicted in figure \ref{fig
ddelta=clasper}.
\begin{figure}[h]
 \centering
 \includegraphics{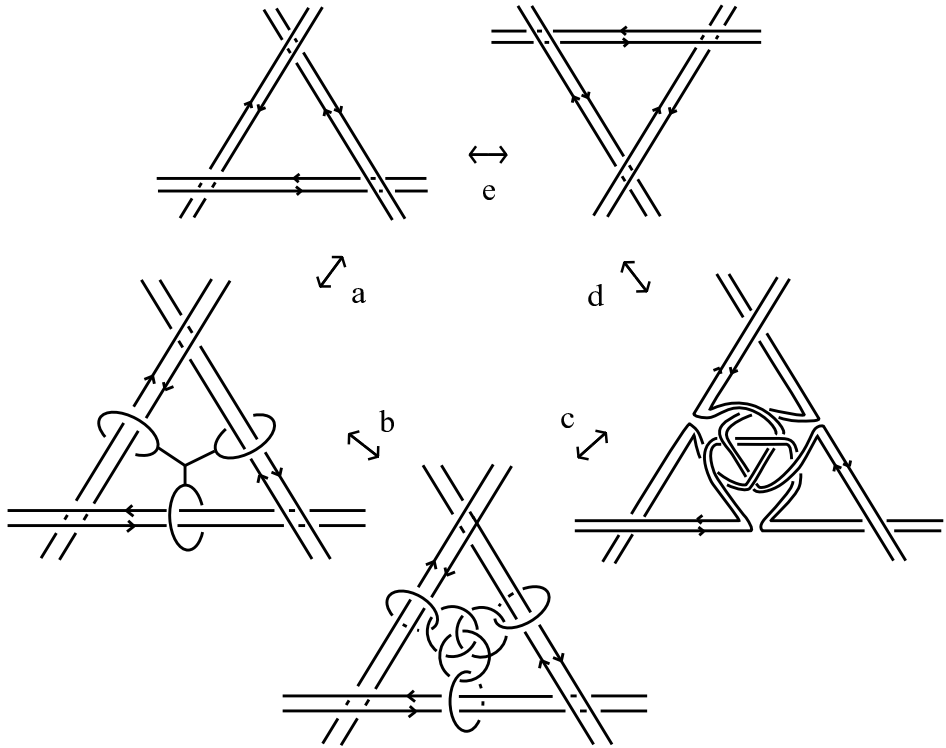}
 \caption{$Y$-Clasper surgery is Doubled-Delta Move}
 \label{fig ddelta=clasper}
\end{figure}
In this figure, the relationship $a$ introduces the desired
$Y$-clasper, $b$ transforms the clasper into its associated link,
$c$ depicts the effect of surgery via handle slide moves, $d$ is
merely the second Reidemeister move, and $e$ is the doubled-delta
move. \
\end{proof}

\subsection{Proof of $ii.\Longrightarrow iii.$}

We actually prove a stronger implication than claimed in the theorem, namely
that the clasper surgery described in $\ ii.$ does not alter the ordered
Seifert matrix. \

\begin{proof}
Consider a neighborhood of the clasp, as in Figure \ref{fig
zoomab}a. \ There are $2k$ strands passing through the leaf, with
$k$ in each direction. \ We can assume the strands' directions
alternate within this neighborhood. \ If they do not alternate, we
can permute the strands by introducing inverse braids just above and
below the $2k$ braid in question.

\begin{figure}[h]
 \centering
 \includegraphics{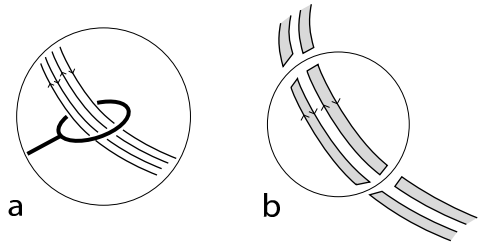}
 \caption{Enlarged View of Clasper Leaf}
 \label{fig zoomab}
\end{figure}

Temporarily cut the strands at the neighborhood's boundary and fill
in the alternating arcs to get bands as in Figure \ref{fig zoomab}b.
Outside the neighborhood, apply Seifert's algorithm to the resulting
link.  Then reattach the new bands at the neighborhood's boundary to
obtain a Seifert surface for the original link where the leaf of the
clasper grabs $k$ bands of the Seifert surface. \

Now clasper surgery is equivalent to tying these bands into
Borromean rings, which doesn't affect the linking of the bands or strands. \
If a basis element of $H_{1}(\Sigma )$ runs through a band, its
pairwise linking with the other basis elements is unchanged; if not,
it is completely unaffected by the surgery. \ Thus the clasper
surgery leaves the Seifert matrix unchanged. $\ $
\end{proof}

\subsection{Proof of $iii.\implies iv.$}

With almost no alteration, the proof of $iii.\implies iv.$ can also be used
to show that two classically S-equivalent knots or links have a Seifert
matrix in common. \

\begin{proof}
If $L$ and $L^{\prime }$ are Strongly S-equivalent then, by definition,
there are Seifert matrices $M$ and $M^{\prime }$ respectively such that $M$
is equivalent to $M^{\prime }$ under a finite sequence of enlargements $%
(\nearrow )$, reductions $(\searrow )$, and congruences $(\cong )$. \ For
example, we can write a sequence of the form
\begin{equation}
M\searrow M_{1}\nearrow M_{2}\cong M_{3}\searrow M_{4}\cong M_{5}\nearrow
M_{6}\searrow M^{\prime }.  \tag{$\ast $}  \label{eqn unordered sequence}
\end{equation}%
In order to prove that $L$ and $L^{\prime }$ have a Seifert matrix in
common, it will be helpful to rewrite the sequence (\ref{eqn unordered
sequence}) so that all the enlargements precede all the reductions, as in
\begin{equation}
M\nearrow \widehat{M_{1}}\cong \widehat{M_{2}}\nearrow \widehat{M_{3}}\cong
\widehat{M_{4}}\nearrow \widehat{M_{5}}\cong \widehat{M_{6}}\searrow
\widehat{M_{7}}\searrow \widehat{M_{8}}\cong M^{\prime }.  \tag{$\ast \ast $}
\label{eqn ordered sequence}
\end{equation}%
The following lemmas establish that such an ordered sequence (\ref{eqn ordered
sequence}) exists for
every pair of strongly S-equivalent matrices and show how the ordered
sequence completes the proof of $iii.\implies iv$.. \

\begin{lemma}
\label{mainlemmaF}$M_{1}\searrow M_{2}\nearrow M_{3}\Longrightarrow
M_{1}\nearrow M_{4}\cong M_{5}\searrow M_{3}$.
\end{lemma}

Before proving this lemma, we first note that the simpler implication $%
M_{1}\searrow M_{2}\nearrow M_{3}\Longrightarrow M_{1}\nearrow M_{4}\searrow
M_{3}$ is not true. \ In the right-hand side, a reduction immediately
follows an enlargement, so the two actions cancel each other and $%
M_{1}=M_{3},$ which is not necessarily true in the left-hand side.

\begin{proof}
Assume that

\begin{equation*}
M_{1}=%
\begin{pmatrix}
&  &  & | & | \\
& V &  & \overrightarrow{y_{1}}^{t} & 0 \\
&  &  & | & | \\
\text{---} & \overrightarrow{x_{1}} & \text{---} & z_{1} & 0 \\
\text{---} & 0 & \text{---} & 1 & 0%
\end{pmatrix}%
,\text{ }M_{2}=V,\text{ and }M_{3}=%
\begin{pmatrix}
&  &  & | & | \\
& V &  & \overrightarrow{y_{3}}^{t} & 0 \\
&  &  & | & | \\
\text{---} & \overrightarrow{x_{3}} & \text{---} & z_{3} & 0 \\
\text{---} & 0 & \text{---} & 1 & 0%
\end{pmatrix}%
.
\end{equation*}%
\ Then let

\begin{equation*}
M_{4}=%
\begin{pmatrix}
&  &  & | & | & | & | \\
& V &  & \overrightarrow{y_{1}}^{t} & 0 & \overrightarrow{y_{3}}^{t} & 0 \\
&  &  & | & | & | & | \\
\text{---} & \overrightarrow{x_{1}} & \text{---} & z_{1} & 0 & 0 & 0 \\
\text{---} & 0 & \text{---} & 1 & 0 & 0 & 0 \\
\text{---} & \overrightarrow{x_{3}} & \text{---} & 0 & 0 & z_{3} & 0 \\
\text{---} & 0 & \text{---} & 0 & 0 & 1 & 0%
\end{pmatrix}%
.
\end{equation*}%
If $M_{4}$ is a $k\times k$ matrix, then let $M_{5}$ be the result of a
change of basis that permutes the $k$th basis element with the $k-2$nd, and
the $k-1$st basis element with the $k-3$rd. \ Then

\begin{equation*}
M_{5}=%
\begin{pmatrix}
&  &  & | & | & | & | \\
& V &  & \overrightarrow{y_{3}}^{t} & 0 & \overrightarrow{y_{1}}^{t} & 0 \\
&  &  & | & | & | & | \\
\text{---} & \overrightarrow{x_{3}} & \text{---} & z_{3} & 0 & 0 & 0 \\
\text{---} & 0 & \text{---} & 1 & 0 & 0 & 0 \\
\text{---} & \overrightarrow{x_{1}} & \text{---} & 0 & 0 & z_{1} & 0 \\
\text{---} & 0 & \text{---} & 0 & 0 & 1 & 0%
\end{pmatrix}%
,
\end{equation*}%
which reduces to $M_{3}$. \
\end{proof}

\begin{lemma}
\label{mainlemmaG}$M_{1}\cong M_{2}\nearrow M_{3}\Longrightarrow
M_{1}\nearrow M_{4}\cong M_{3}$.
\end{lemma}

\begin{proofnosquare}
If $M_{2}=P^{t}M_{1}P$, and%
\begin{eqnarray*}
M_{3} &=&%
\begin{pmatrix}
&  &  & | & | \\
& P^{t}M_{1}P &  & \overrightarrow{y_{3}}^{t} & 0 \\
&  &  & | & | \\
\text{---} & \overrightarrow{x_{3}} & \text{---} & z_{3} & 0 \\
\text{---} & 0 & \text{---} & 1 & 0%
\end{pmatrix}%
\text{, } \\
\text{let \ }M_{4} &=&%
\begin{pmatrix}
&  &  & | & | \\
& M_{1} &  & (P^{t})^{-1}\cdot \overrightarrow{y_{3}}^{t} & 0 \\
&  &  & | & | \\
\text{---} & \overrightarrow{x_{3}}\cdot P^{-1} & \text{---} & z_{3} & 0 \\
\text{---} & 0 & \text{---} & 1 & 0%
\end{pmatrix}%
\end{eqnarray*}%
be an enlargement of $M_{1}$. \ Now let $Q=\left(
\begin{smallmatrix}
P & 0 \\
0 & I%
\end{smallmatrix}%
\right) $. \ Then%
\begin{eqnarray*}
Q^{t}M_{4}Q &=&%
\begin{pmatrix}
P^{t} & 0 \\
0 & I%
\end{pmatrix}%
\begin{pmatrix}
&  &  & | & | \\
& M_{1} &  & (P^{t})^{-1}\cdot \overrightarrow{y_{3}}^{t} & 0 \\
&  &  & | & | \\
\text{---} & \overrightarrow{x_{3}}\cdot P^{-1} & \text{---} & z_{3} & 0 \\
\text{---} & 0 & \text{---} & 1 & 0%
\end{pmatrix}%
\begin{pmatrix}
P & 0 \\
0 & I%
\end{pmatrix}
\\
&=&%
\begin{pmatrix}
&  &  & | & | \\
& P^{t}M_{1} &  & P^{t}\cdot \lbrack (P^{t})^{-1}\cdot \overrightarrow{y_{3}}%
^{t}] & 0 \\
&  &  & | & | \\
\text{---} & \overrightarrow{x_{3}}\cdot P^{-1} & \text{---} & z_{3} & 0 \\
\text{---} & 0 & \text{---} & 1 & 0%
\end{pmatrix}%
\begin{pmatrix}
P & 0 \\
0 & I%
\end{pmatrix}
\\
&=&%
\begin{pmatrix}
&  &  & | & | \\
& P^{t}M_{1}P &  & P^{t}\cdot \lbrack (P^{t})^{-1}\cdot \overrightarrow{y_{3}%
}^{t}] & 0 \\
&  &  & | & | \\
\text{---} & [\overrightarrow{x_{3}}\cdot P^{-1}]\cdot P & \text{---} & z_{3}
& 0 \\
\text{---} & 0 & \text{---} & 1 & 0%
\end{pmatrix}
\\
&=&%
\begin{pmatrix}
&  &  & | & | \\
& P^{t}M_{1}P &  & \overrightarrow{y_{3}}^{t} & 0 \\
&  &  & | & | \\
\text{---} & \overrightarrow{x_{3}} & \text{---} & z_{3} & 0 \\
\text{---} & 0 & \text{---} & 1 & 0%
\end{pmatrix}%
=M_{3}
\end{eqnarray*}%
\hfill $\square $
\end{proofnosquare}

\begin{lemma}
\label{mainlemmaA1}Any sequence of relations between Strongly S-equivalent
matrices can be rewritten so that all enlargements come before all
reductions, as in the sequence (\ref{eqn ordered sequence}).
\end{lemma}

\begin{proof}
The proof of Lemma \ref{mainlemmaA1} follows from an induction argument
using Lemmas \ref{mainlemmaF} and \ref{mainlemmaG}. \ For simplicity, we
will use strings of the symbols \{$\nearrow ,\searrow ,\cong $\} to denote
enlargements, reductions, and congruences, respectively, while omitting
explicit reference to the matrices.

\begin{itemize}
\item Base Case: \ $\searrow $ $\nearrow $ $\implies $ $\nearrow $ $\cong $ $%
\searrow $ (Lemma \ref{mainlemmaF}).

\item Inductive Step: \ Find the first enlargement that is preceded by a reduction.
 \ Note that the sequence preceding this enlargement is
arranged as desired. \ There are two cases:

\begin{enumerate}
\item If this enlargement is immediately preceded by a congruence, apply
Lemma \ref{mainlemmaG} to replace $\cong $ $\nearrow $ with $\nearrow $ $%
\cong $. \ Repeat Lemma \ref{mainlemmaG} until the enlargement is
immediately preceded by a reduction.

\item If this enlargement is immediately preceded by a reduction, then by
Lemma \ref{mainlemmaF}, $\searrow $ $\nearrow $ can be replaced with $%
\nearrow $ $\cong $ $\searrow $. \
\end{enumerate}

Continue the two steps above until the enlargement is immediately preceded
by another enlargement. \
\end{itemize}

We have reduced the number of out-of-order \{$\nearrow ,\searrow ,\cong $\}
by one without increasing the number of enlargements or reductions. \
Continue the inductive step until there is no enlargement preceded by a
reduction. \
\end{proof}

We can now finish the proof of $iii.\implies iv.$. \ If $L$ and $L^{\prime }$
are Strongly S-equivalent, then by Lemma \ref{mainlemmaA1} we may assume
there exists a sequence of relations $(\nearrow ,\searrow ,\cong )$ between $%
M$ and $M^{\prime }$ where all enlargements precede all reductions. \ Note
that:

\begin{itemize}
\item If $M$ is an ordered Seifert matrix for the link $L$ with respect
to Seifert surface $\Sigma $ and basis $\beta $, and if $M\nearrow M^{\prime
}$, then $M^{\prime }$ is also an ordered Seifert matrix for $L$ with
respect to a 1-handle enlargement $\widehat{\Sigma }$ of $\Sigma $ and the
corresponding new basis $\beta \cup \{a_{1},a_{2}\}$. \ This follows from
Proposition \ref{prop enlargement}.

\item If $M$ is an ordered Seifert matrix for the link $L$, and $M\cong
M^{\prime }$, then $M^{\prime }$ is also an ordered Seifert matrix for $L$,
with respect to the same Seifert surface and a new basis as prescribed by
the congruence.
\end{itemize}

Using the two facts above, we can \textquotedblleft work inwards\textquotedblright\  
from both ends of the ordered sequence to show that $L$
and $L^{\prime }$ have a common Seifert matrix (though not necessarily $M$
or $M^{\prime }$) for some Seifert surfaces and bases. \ Starting with $M$
and working from left to right, each enlargement or congruence yields a new
Seifert matrix for $L$. \ This terminates when we reach the first reduction.
\ Similarly, since $M_{i}\searrow M^{\prime }\iff $ $M^{\prime }\nearrow
M_{i}$, starting with $M^{\prime }$ and working from right to left, each
\emph{reduction} or congruence yields a new Seifert matrix for $L^{\prime }$%
. \ In the ordered sequence labelled (\ref{eqn ordered sequence}) above, $%
\widehat{M_{5}}$ is a common Seifert matrix for $L$ and $L^{\prime }$, as is $%
\widehat{M_{6}}$. \
\end{proof}

\subsection{Proof of $iv.\implies i.$}

The following proposition contains the bulk of content of the Main Theorem.
\ This is the primary step that distinguishes the proof of $iv.\implies i.$
from Naik-Stanford's proof and uses the extra hypotheses of Strong
S-equivalence:

\begin{proposition}
\label{prop nice disk-band}If two $m$-component links $L_{0}$ and $L_{1}$
have the same ordered Seifert matrix $M$ with respect to Seifert surfaces $%
\Sigma _{0}$ and $\Sigma _{1}$ and ordered bases $\beta _{0}$ and $\beta
_{1} $ of $H_{1}(\Sigma _{0})$ and $H_{1}(\Sigma _{1})$, respectively, then
it is possible to arrange $\Sigma _{0}$ and $\Sigma _{1}$ into disk-band
form \emph{and} to find new semi-symplectic\ bases $\gamma _{0}$ and $\gamma
_{1}$ for $\Sigma _{0}$ and $\Sigma _{1}$ that give rise to a new shared
ordered Seifert matrix $N$ for both $L_{0}$ and $L_{1}$. \
\end{proposition}

\begin{proof}[Proof of Proposition]
We start with surfaces $\Sigma _{0}$ and $\Sigma _{1}$ and bases $\beta _{0}$
and $\beta _{1}$ of $H_{1}(\Sigma _{0})$ and $H_{1}(\Sigma _{1})$,
respectively. \ In transforming $\Sigma _{0}$ and $\Sigma _{1}$ into
disk-band form, it is important that we keep track of the new semi-symplectic 
bases $\gamma _{0}$ and $\gamma _{1}$ in terms
$\beta _{0}$ and $\beta _{1}$. \ This means understanding the homeomorphisms
involved in the transformation. \

Both $\Sigma _{0}$ and $\Sigma _{1}$ have $m$ boundary components, and since
they share an ordered Seifert matrix, both have the same genus, say $g$. \ Let $F_{g}$
be an abstract surface of genus $g$ with $m$ boundary components,
specifically realized as a disk with bands as in Figure \ref{fig diskband} (%
\emph{i.e.} the string link from Figure \ref{fig diskband stringlink w/basis}
is trivial),
\begin{figure}[h]
 \centering
 \includegraphics{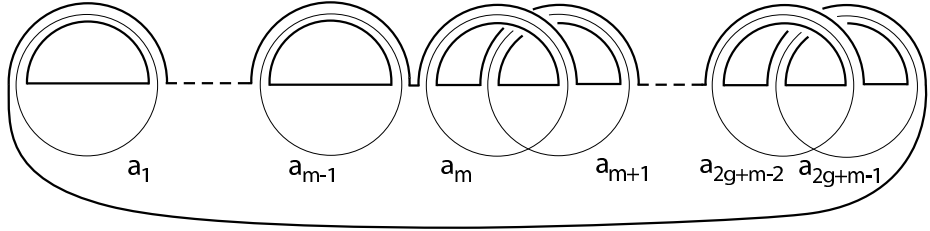}
 \caption{The surface $F_{g}$}
 \label{fig diskband}
\end{figure}
where the $\{a_{i}\}$ form an ordered basis for $H_{1}(F_{g})$. \
The intersection form for $\{a_{i}\}$ is represented by the block
matrix $X=\left(
\begin{smallmatrix}
0 & 0 \\
0 & Sym%
\end{smallmatrix}%
\right) $, where $Sym=\bigoplus\limits_{j=1}^{g}\left(
\begin{smallmatrix}
0 & -1 \\
1 & 0%
\end{smallmatrix}%
\right) $. \

Choose orientation-preserving homeomorphisms $\Phi _{i}:F_{g}\rightarrow
\Sigma _{i}$ from the pure mapping class group, \emph{i.e.} the boundary
components are fixed, pointwise. \ Assume that the boundary component of $%
F_{g}$ that is parallel to $a_{j}$ is sent by $\Phi _{i}$ to the $j$th
component of $L_{i}$. \ Then
\begin{eqnarray*}
\Phi _{i}(\{a_{i}\}) &=&\{\Phi _{i}(a_{1}),\ldots ,\Phi _{i}(a_{m-1+2g})\} \\
&=&\{\beta _{i,1},\beta _{i,2},\ldots ,\beta _{i,m-1},\Phi
_{i}(a_{m}),\ldots ,\Phi _{i}(a_{m-1+2g})\}
\end{eqnarray*}%
is also a basis for $H_{1}(\Sigma _{i})$. \ So there exist invertible
matrices $A_{i}=$ $\left(
\begin{smallmatrix}
I & \ast \\
0 & \ast%
\end{smallmatrix}%
\right) $ such that $N_{i}=A_{i}^{t}MA_{i}$ denote the ordered Seifert matrices of $%
\Sigma _{i}$ with respect to the new bases $\Phi _{i}(\{a_{i}\})$. \ These
homeomorphisms describe an explicit way to put $\Sigma _{0}$ and $\Sigma
_{1} $ into disk band form as in Figure \ref{fig diskband stringlink w/basis}%
, where the two surfaces differ only by string links $\Lambda _{0}$, $%
\Lambda _{1}$ of bands of the Seifert surfaces. \ The new ordered Seifert matrices $%
N_{0}$ and $N_{1}$ correspond to these \textquotedblleft
new\textquotedblright\ surfaces, but $N_{0} \neq N_{1}$.  In the construction process we lost
the critical hypothesis that the ordered Seifert matrices for $L_{0}$ and $%
L_{1}$ were equal. \ We need to show that $N_{0}$ is also an ordered Seifert matrix
for $\Sigma _{1}$, with respect to a different basis. \

By definition, if we have a Seifert matrix $N$ determined by a given basis,
then $N-N^{t}$ is the intersection form for that same basis. \ Since
homeomorphisms of surfaces preserve the intersection properties of their
basis elements, $\Phi _{i}(\{a_{i}\})$ will have the same intersection
properties as $\{a_{i}\}$, and thus we have that $N_{i}-N_{i}^{t}=X=\left(
\begin{smallmatrix}
0 & 0 \\
0 & Sym%
\end{smallmatrix}%
\right) $. \ Using this fact, we construct a matrix $C=A_{1}^{-1}A_{0}$ that
we show stabilizes $X$, in that $X=C^{t}XC$. \ First we show that $%
N_{0}=C^{t}N_{1}C$:%
\begin{eqnarray*}
N_{0} &=&A_{0}^{t}MA_{0} \\
&=&A_{0}^{t}((A_{1}^{t})^{-1}MA_{1}^{-1})A_{0} \\
&=&(A_{1}^{-1}A_{0})^{t}N_{1}(A_{1}^{-1}A_{0}) \\
&=&C^{t}N_{1}C
\end{eqnarray*}%
Then,%
\begin{eqnarray*}
X &=&N_{0}-N_{0}^{t} \\
&=&C^{t}N_{1}C-C^{t}N_{1}^{t}C \\
&=&C^{t}(N_{1}-N_{1}^{t})C \\
&=&C^{t}XC.
\end{eqnarray*}

Now the key step is to find a homeomorphism of the pure
mapping class group that induces $C$ with respect to the ordered basis $\Phi
_{1}(\{a_{i}\})$. This, in turn, will prove that $N_{0}$ is an ordered Seifert 
matrix for $\Sigma _{1}$.  (Actually the pure mapping class group, which fixes the
boundary components pointwise, is stronger than necessary. \ Fixing the
boundary component-wise would be sufficient; however, the ordinary mapping
class group only fixes the boundary set-wise.) \ For the case of knots, this
key step is easy. \ There $C$ is symplectic, and thus is well-known to be induced
by such a homeomorphism, since the map from the pure mapping class group to
the symplectic group is surjective \cite[pp. 178, 355-6]{Magnus}. \ For
links however, $C$ is not symplectic, and it takes several steps to find a 
homeomorphism inducing $C$.

\begin{lemma}
\label{mainlemmaC}The matrix $C$ is of the form $\left(
\begin{smallmatrix}
I & B \\
0 & S%
\end{smallmatrix}%
\right) $ where $S$ is a symplectic matrix.
\end{lemma}

\begin{proof}
By construction, $C=A_{1}^{-1}A_{0}$. \ We know by definition that $A_{0}$
and $A_{1}$ are of the form $A_{i}=$ $\left(
\begin{smallmatrix}
I & \ast \\
0 & \ast%
\end{smallmatrix}%
\right) $. \

First we want to show that $A_{1}^{-1}$ is also of the form $\left(
\begin{smallmatrix}
I & \ast \\
0 & \ast%
\end{smallmatrix}%
\right) $. \ Let $A_{1}=$ $\left(
\begin{smallmatrix}
I & Y \\
0 & Z%
\end{smallmatrix}%
\right) $. \ Suppose $\left(
\begin{smallmatrix}
P & Q \\
R & T%
\end{smallmatrix}%
\right) $ is an inverse for $A_{1}$. \ Then

\begin{equation*}
\begin{pmatrix}
I & 0 \\
0 & I%
\end{pmatrix}%
=%
\begin{pmatrix}
P & Q \\
R & T%
\end{pmatrix}%
\begin{pmatrix}
I & Y \\
0 & Z%
\end{pmatrix}%
=%
\begin{pmatrix}
P & PY+QZ \\
R & RY+TZ%
\end{pmatrix}%
,
\end{equation*}%
so $R=0$ and $P=I$. \ Thus $A_{1}^{-1}$ must be of the form $\left(
\begin{smallmatrix}
I & \ast \\
0 & \ast%
\end{smallmatrix}%
\right) $. \

From this it is easy to see that

\begin{equation*}
C=A_{1}^{-1}A_{0}=%
\begin{pmatrix}
I & Q \\
0 & T%
\end{pmatrix}%
\begin{pmatrix}
I & V \\
0 & W%
\end{pmatrix}%
=%
\begin{pmatrix}
I & V+QW \\
0 & TW%
\end{pmatrix}%
\end{equation*}%
is also of the form $\left(
\begin{smallmatrix}
I & \ast \\
0 & \ast%
\end{smallmatrix}%
\right) $.

\ Now we need to demonstrate that the lower right block of $C$ is
symplectic. \ For this, we let $C=\left(
\begin{smallmatrix}
I & B \\
0 & S%
\end{smallmatrix}%
\right) $. \ Then

\begin{eqnarray*}
C^{t}XC &=&%
\begin{pmatrix}
I & 0 \\
B^{t} & S^{t}%
\end{pmatrix}%
\begin{pmatrix}
0 & 0 \\
0 & Sym%
\end{pmatrix}%
\begin{pmatrix}
I & B \\
0 & S%
\end{pmatrix}
\\
&=&%
\begin{pmatrix}
0 & 0 \\
0 & S^{t}\cdot Sym%
\end{pmatrix}%
\begin{pmatrix}
I & B \\
0 & S%
\end{pmatrix}
\\
&=&%
\begin{pmatrix}
0 & 0 \\
0 & S^{t}\cdot Sym\cdot S%
\end{pmatrix}%
=%
\begin{pmatrix}
0 & 0 \\
0 & Sym%
\end{pmatrix}%
=X.
\end{eqnarray*}

Recall that $Sym=\bigoplus\limits_{j=1}^{g}\left(
\begin{smallmatrix}
0 & -1 \\
1 & 0%
\end{smallmatrix}%
\right) $. Since $S^{t}\cdot Sym\cdot S=Sym$, $S$ is a symplectic matrix.
\end{proof}

\begin{lemma}
\label{mainlemmaD}The matrix $D=\left(
\begin{smallmatrix}
I & 0 \\
0 & S%
\end{smallmatrix}%
\right) $, relative to any fixed semi-symplectic basis, is induced by a
homeomorphism of the pure mapping class group.
\end{lemma}

\begin{proof}
We want to consider $D$ to represent an action on $H_{1}(\Sigma _{1})$. \ If
$\Sigma $ were a once-punctured surface (in other words, if our link were a
knot), then $D$ would be symplectic and this action would then be induced by
a homeomorphism from the pure mapping class group \cite{Magnus}. \

Our Seifert surface $\Sigma _{1}$ is an $m$-punctured surface. \ Choose a
representative curve $\overline{\beta _{i}}\subset \Sigma _{1}$ for each of
the last $2g$ basis elements of $\beta _{1}$, where $m\leq i\leq 2g+m-1$,
such that the geometric intersections among representatives respect the
algebraic intersections determined by the ordered Seifert matrix. \

Temporarily cap off the $m$ boundary components with disks $D_{1},\ldots
,D_{m}$, forming a new closed surface $\overline{\Sigma }$. \ Consider a
larger disk $U\subset \overline{\Sigma }$ that encompasses all the
smaller disks $D_{1},\ldots ,D_{m}$, but does not intersect any of 
$\overline{\beta }_{m},\ldots ,\overline{\beta }_{2g+m-1}$. \ This is
possible. \ Since the $2g$ basis elements are semi-symplectic, we can cut
$\Sigma $ along these curves to obtain an $m$-punctured $2g$-gon. \ The disk $U$ can
be chosen in the interior of this $2g$-gon such that it will encompass all of
the disks $D_{i}$ (a pushoff of the $2g$-gon's boundary will suffice). \

Now $\overline{\Sigma }-U$, a subset of $\Sigma _{1}$, is a once-punctured surface
with symplectic basis $\overline{\beta }_{m},\ldots ,\overline{\beta }%
_{2g+m-1}$. \ Where $D=\left(
\begin{smallmatrix}
I & 0 \\
0 & S%
\end{smallmatrix}%
\right) $ represents an action on $H_{1}(\Sigma _{1})$, $S$ represents an
action on $H_{1}(\overline{\Sigma }-U)$ and corresponds to a homeomorphism $%
\overline{g}$ of the pure mapping class group by the surjective map
mentioned above. \

We can extend $\overline{g}$ to $g$, which agrees with $\overline{g}$ on $%
\overline{\Sigma }-U\subset \Sigma _{1}$ and fixes the remaining $%
U-\bigcup\limits_{j=1}^{m}D_{j}\subset \Sigma _{1}$. \ We now have that the
homeomorphism $g$ induces the action $D=\left(
\begin{smallmatrix}
I & 0 \\
0 & S%
\end{smallmatrix}%
\right) $ on $H_{1}(\Sigma _{1})$.\
\end{proof}

\begin{lemma}
\label{mainlemmaE}There is a matrix $E$ such that $C=DE$, and $E$ can be
taken to be a product of elementary matrices $E_{i,j}$ where each $E_{i,j}$
is induced by a homeomorphism of the pure mapping class group.
\end{lemma}

\begin{proof}
Observe that $\left(
\begin{smallmatrix}
I & B \\
0 & S%
\end{smallmatrix}%
\right) = \left(
\begin{smallmatrix}
I & 0 \\
0 & S%
\end{smallmatrix}%
\right) \left(
\begin{smallmatrix}
I & B \\
0 & I%
\end{smallmatrix}%
\right) $, or $C = D \left(
\begin{smallmatrix}
I & B \\
0 & I%
\end{smallmatrix}%
\right) $.  Let $E=\left(
\begin{smallmatrix}
I & B \\
0 & I%
\end{smallmatrix}%
\right) $. \ Right-multiplication by this matrix $E$ can be achieved by a
product of several elementary matrices $E_{i,j}$, which we will proceed to
define. \

Note that $B$ is an $(m-1)\times 2g$ block.  To better illustrate the 
construction of the $E_{i,j}$, we first let $B_{i,j}$ be the $%
(m-1)\times 2g$ matrix that has the $b_{i,j}$ entry of $B$ in the $i,j$th
spot and zeros elsewhere. \ One can verify that $E=\left(
\begin{smallmatrix}
I & B \\
0 & I%
\end{smallmatrix}%
\right) =\prod\limits_{i,j=1}^{m-1}\left(
\begin{smallmatrix}
I & B_{i,j} \\
0 & I%
\end{smallmatrix}%
\right) $. \ Define an elementary $(2g+m-1)\times (2g+m-1)$ matrix $E_{i,j}$
to be the identity matrix with $+1$ in the $i,j+m-1$ spot, where $1\leq
i\leq m-1$, and $1\leq j\leq 2g$. \ For example, if $m=4$ and $g=2$ then

\begin{equation*}
E_{1,2}=%
\begin{pmatrix}
I &
\begin{array}{cccc}
0 & 1 & 0 & 0 \\
0 & 0 & 0 & 0 \\
0 & 0 & 0 & 0%
\end{array}
\\
0 & I%
\end{pmatrix}%
=%
\begin{pmatrix}
1 & 0 & 0 & 0 & 1 & 0 & 0 \\
0 & 1 & 0 & 0 & 0 & 0 & 0 \\
0 & 0 & 1 & 0 & 0 & 0 & 0 \\
0 & 0 & 0 & 1 & 0 & 0 & 0 \\
0 & 0 & 0 & 0 & 1 & 0 & 0 \\
0 & 0 & 0 & 0 & 0 & 1 & 0 \\
0 & 0 & 0 & 0 & 0 & 0 & 1%
\end{pmatrix}%
\end{equation*}

It should also be noted that for any integer $k$,

\begin{equation*}
(E_{1,2})^{k}=%
\begin{pmatrix}
I &
\begin{array}{cccc}
0 & k & 0 & 0 \\
0 & 0 & 0 & 0 \\
0 & 0 & 0 & 0%
\end{array}
\\
0 & I%
\end{pmatrix}%
\end{equation*}

Then $C=DE$, where $E=\prod\limits_{i=1}^{m-1}\prod%
\limits_{j=1}^{2g}(E_{i,j})^{b_{i,j}}$. \

To show that these elementary matrices $E_{i,j}$ are induced by
homeomorphisms of the pure mapping class group, note that each $E_{i,j}$
fixes the first $m-1$ basis elements and sends $b_{m-1+j}$ to $%
b_{m-1+j}+b_{i}$. \ An example of this is shown below for $E_{1,3}$ where $%
m=4$ and $g=2$:

\begin{eqnarray*}
\begin{pmatrix}
I &
\begin{array}{cccc}
0 & 0 & 1 & 0 \\
0 & 0 & 0 & 0 \\
0 & 0 & 0 & 0%
\end{array}
\\
0 & I%
\end{pmatrix}%
\begin{pmatrix}
0 \\
1 \\
0 \\
0 \\
0 \\
0 \\
0%
\end{pmatrix}
&=&%
\begin{pmatrix}
0 \\
1 \\
0 \\
0 \\
0 \\
0 \\
0%
\end{pmatrix}
\\
\text{and \ \ \ }%
\begin{pmatrix}
I &
\begin{array}{cccc}
0 & 0 & 1 & 0 \\
0 & 0 & 0 & 0 \\
0 & 0 & 0 & 0%
\end{array}
\\
0 & I%
\end{pmatrix}%
\begin{pmatrix}
0 \\
0 \\
0 \\
0 \\
0 \\
1 \\
0%
\end{pmatrix}
&=&%
\begin{pmatrix}
1 \\
0 \\
0 \\
0 \\
0 \\
1 \\
0%
\end{pmatrix}%
\text{.}
\end{eqnarray*}%
We need to find a homeomorphism that induces $E_{i,j}$, that is, one that
exactly takes $b_{k}$ to $b_{k}+b_{i}$ for $k=j+m-1$ (and $1\leq i\leq
m-1,1\leq j\leq 2g$), while fixing all the other basis elements. \ These
homeomorphisms can be realized as simple Dehn twists, and are illustrated in
figure \ref{fig dehn}.

\begin{figure}[h]
 \centering
 \includegraphics{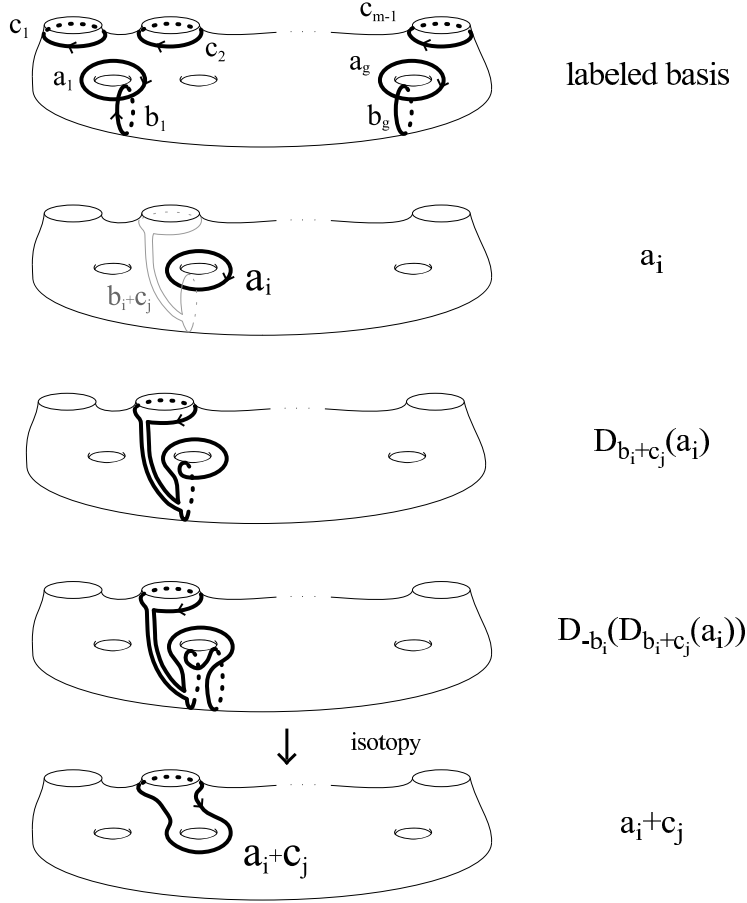}
 \caption{Dehn twists $D_{-b_{i}}(D_{b_{i}+c_{j}}(a_{i}))=a_{i}+c_{j}$.}
 \label{fig dehn}
\end{figure}

In each of these figures, the surface $\Sigma _{1}$ is represented
schematically as a genus $2g$ surface with $m$ boundary components.
\ Pairs of the standard symplectic basis elements are labelled
$\{a_{i},b_{i}\}$, while the boundary basis elements are labelled
$c_{i}.$ \ The notation $D_{b}(a)$ is used to represent the effect
on the curve $a$ of performing a Dehn twist about the curve $b$. \
Using this notation, we show that $a_{i}$ can be taken to
$a_{i}+c_{j}$ by
the composition of Dehn twists $%
D_{-b_{i}}(D_{b_{i}+c_{j}}(a_{i}))=a_{i}+c_{j}$. \ Similarly, $%
D_{-a_{i}}(D_{a_{i}+c_{j}}(b_{i}))=b_{i}+c_{j}$. \

Since each $E_{i,j}$ fixes the first $m-1$ basis elements, $E_{i,j}$ is
induced by a homeomorphism $f_{i,j}$ that fixes the first $m-1$ boundary
components of $\Sigma _{i}$ and hence fixes all $m$ boundary components. \
All that is required for the sake of $C$ is that the boundary be fixed,
component-wise. \ However, as the $f_{i,j}$ are constructed from these
simple Dehn twists, they can easily be taken to fix the boundary components
point-wise. \ Now, $E$, as a product of elementary matrices $E_{i,j}$, is
induced by the corresponding composition of the $f_{i,j}$. \ That is, the
action of $E$ on $H_{1}(\Sigma _{1})$ is induced by the composition $f=%
\underset{i=1}{\overset{m-1}{\circ }}\overset{2g}{\underset{j=1}{\circ }}%
\overset{b_{i,j}}{\underset{k=1}{\circ }}f_{i,j}$ of homeomorphisms $f_{i,j}$%
. \
\end{proof}

Together, Lemmas \ref{mainlemmaD} and \ref{mainlemmaE} imply that $C$ is
induced by a homeomorphism $h=f\circ g$ of the pure mapping class group. 
Now we can use the homeomorphisms $\Phi _{0}:F_{g}\rightarrow \Sigma _{0}$
and $h\circ \Phi _{1}:F_{g}\rightarrow \Sigma _{1}$ to put $\Sigma _{0}$ and
$\Sigma _{1}$ into disk and band form as in Figure \ref{fig diskband
stringlink}. \ The basis elements $h\circ \Phi _{1}(\{a_{i}\})$ are the
columns of the matrix $C$. \ Where $N_{1}$ was the Seifert matrix for $%
\Sigma _{1}$ with respect to the basis $\Phi _{1}(\{a_{i}\})$, now $%
N_{0}=h_{\ast }^{t}N_{1}h_{\ast }=C^{t}N_{1}C$ is the Seifert matrix for $%
\Sigma _{1}$ with respect to the basis $h\circ \Phi _{1}(\{a_{i}\})$. \ Thus
the Seifert surfaces $\Sigma _{0}$ and $\Sigma _{1}$, together with the
semi-symplectic\ bases $\Phi _{0}(\{a_{i}\})$ and $h\circ \Phi
_{1}(\{a_{i}\})$, respectively, give rise to a shared ordered Seifert matrix
$N_{0}$ for $L_{0}$ and $L_{1}$. \
\end{proof}

Finally, the tools are in place to prove the final implication of the Main
Theorem.

\begin{proof}[Proof of $iv.\implies i.$]
Suppose $L_{0}$ and $L_{1}$ have the same ordered Seifert matrix $M$ with
respect to Seifert surfaces $\Sigma _{0}$ and $\Sigma _{1}$ and ordered
bases $\beta _{0}$ and $\beta _{1}$ of $H_{1}(\Sigma _{0})$ and $%
H_{1}(\Sigma _{1})$, respectively. \ Proposition \ref{prop nice disk-band}
states that the Seifert surfaces $\Sigma _{0}$ and $\Sigma _{1}$ can be
arranged in disk-band form as in Figure \ref{fig diskband stringlink}, with
new semi-symplectic\ bases $\gamma _{0}$ and $\gamma _{1}$ that each give
rise to a shared ordered Seifert matrix $N$ for $L_{0}$ and $L_{1}$.

\begin{figure}[h]
 \centering
 \includegraphics{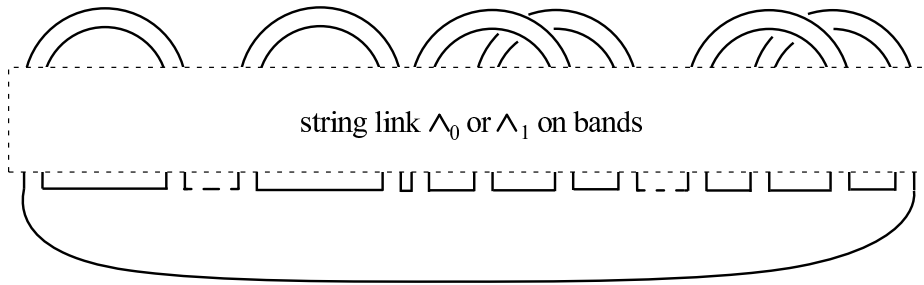}
 \caption{$\Sigma _{0}$ and $\Sigma _{1}$ differ
by string links $\Lambda _{0}$ and $\Lambda _{1}$ on their bands}
 \label{fig diskband stringlink}
\end{figure}

In their new disk-band form, the Seifert surfaces $\Sigma _{0}$ and $\Sigma
_{1}$ differ only by the string links $\Lambda _{0}$ and $\Lambda _{1}$ of
bands. \ Because $\gamma _{0}$ and $\gamma _{1}$ were chosen to be
semi-symplectic, the new basis elements run straight through each band. \
The sets of pairwise linking numbers for $\Lambda _{0}$ and $\Lambda _{1}$
are equal, both being determined by the matrix $N$. \ By
Murakami-Nakanishi's theorem \cite{Murakami}, $\Lambda _{0}$ and $\Lambda
_{1}$ are related by a sequence of (single) delta moves. \ Therefore our
original links $L_{0}$ and $L_{1}$ are related by a sequence of
doubled-delta moves, having two oppositely oriented strands for each band of
$\Lambda _{0}$ or $\Lambda _{1}$. \

Although the issue of framing on the bands of $\Sigma _{0}$ and $\Sigma _{1}$
has not yet been addressed, it is clear that the delta move doesn't change the
framing on any strand of a string link, and the doubled-delta move doesn't
alter the self-linking of any of the bands. \ Moreover, the framing of each
band corresponds to the self-linking of the basis element running through
that band and is thus the same for $\Sigma _{0}$ as for $\Sigma _{1}$, both
being determined by the diagonal entries of $N$. \
\end{proof}

\end{document}